\documentclass[journal, letterpaper]{IEEEtran}

\usepackage{graphicx}
\usepackage{amsmath,amstext,amssymb,amsopn,amsthm}
\newtheorem{definition}{Definition}

\usepackage{url,verbatim}
\usepackage{mathtools}
\usepackage{enumerate}
\usepackage{float}
\usepackage{url}        
\usepackage{textgreek}
\usepackage{listings}
\usepackage{csvsimple}
\usepackage{longtable}
\usepackage{caption}

\usepackage{eucal,mathrsfs,dsfont}

\newtheorem{theorem}{Theorem}[section]

\theoremstyle{definition}

\usepackage{algorithm}
\usepackage[noend]{algpseudocode}

\numberwithin{equation}{section}

\usepackage{color,graphicx,tikz}
\usepackage{subfig}

\begin{document}

	\title{Pinned Billiard Balls Simulation\\ (WXML Autumn 2023 report)}
	\author{Bella Cochran, Anish Gupta, Colin Holzman-Klima, Wanchaloem Wunkaew\\ Advised by Krzysztof Burdzy and Raghavendra Tripathi\\\textit{University of Washington. Seattle, WA, United States.}}
	\maketitle

\begin{abstract}
Systems of pinned billiard balls serve as simplified models of collisions, where all particles remain fixed in their positions while their (pseudo-)velocities evolve in accordance with the laws of conservation of energy and momentum. For some families of ball configurations, Athreya, Burdzy, and Duarte have established the maximum upper bound for the number of pseudo-collisions, thereby demonstrating that the number of collisions is finite (\cite{AthreyaBurdzyDuarte2018PinnedBilliardFoldings}). The result has been extended to all ball configurations in \cite{BuDu2022}. In this project, we do extensive simulations to study two specific configurations. First, we consider balls arranged in a half-space and assign a single ball an inward (pseudo-) velocity. Simulations suggest that in the long run, most of the energy is concentrated near the boundary. Second, when the balls are arranged on a flat torus, we find that in the stationary regime, the distributions of the velocity components are i.i.d. normal. Additionally, we find that the components of the velocities in the direction of impact between two touching balls are uncorrelated.
\end{abstract}

\section{Introduction}
\label{Sec:Intro}
The dynamics of billiard balls, governed by the principles of classical mechanics, have long fascinated both mathematicians and physicists alike \cite{Tabachnikov2005GeometryAndBilliards} \cite{KozlovAndTreshchev1991BilliardsGenetic}, with sub-fields such as the study of semi-dispersing billiards having been investigated with geometric approaches \cite{BuragoFerlegerKononenko1998GeometricApproachSemiDispersingBilliards}, in pursuit of finding the global bounds on collisions of said system \cite{BuragoFerlegerKononenko1998GlobalBoundNumCollisionsSemiDispersingBilliards}, and in the context of both Riemannian manifolds \cite{BuragoFerlegerKononenko1998EstimateOnNumCollisionsSemiDispersingBilliards} and Banach spaces \cite{Goebel2018BallsInBanachSpaces}. Traditionally, these studies have focused on the behavior of moving balls, as in \cite{Galperin1981MovingLocalRepellingParticles} and \cite{Galperin1983OnMovingLocalRepellingParticles}, under the influence of elastic collisions within defined boundaries. However, the exploration of billiard ball dynamics takes an interesting turn with the introduction of the ``pinned balls" model. This novel concept, which centers on balls that are fixed in position but possess pseudo-velocities, offers a unique perspective on collision dynamics. 

In this pinned model, the balls remain stationary, but their pseudo-velocities evolve following the laws of conservation of energy and momentum, akin to those governing totally elastic collisions in moving balls. This approach simplifies the complex interaction of colliding particles making it an ideal framework for theoretical exploration and computational simulation. Thus the study of pinned billiard balls serves as a simplified but potent model for understanding the fundamental aspects of collision dynamics in constrained systems. Our project delves into the dynamics of pinned billiard balls, building upon the foundational works that have established the theoretical underpinnings of this model \cite{AthreyaBurdzyDuarte2018PinnedBilliardFoldings}. We focus on extensive simulations to study two distinct environments: firstly, balls arranged in a half-space configuration, and secondly, balls organized on a flat torus. These simulations aim to investigate the equilibrium states and energy distribution of these setups, providing insights into the behavior of the system under different boundary conditions. In the half-space configuration, our simulations suggest a concentration of energy near the boundary, highlighting the influence of spatial constraints on the system's dynamics. Separately, in the flat torus arrangement, we observe that the velocity components tend to distribute normally in the stationary regime, reminiscent of Markov processes, \cite{BurdzyAndWhite2008MarkovWithProductFormStationaryDistributions}. Further, our findings indicate that the components of the velocities in the direction of impact between touching balls are uncorrelated. Through our study, we aim to contribute to the broader understanding of collision dynamics in constrained environments, offering new perspectives on the behavior of pseudo-collisions and energy distribution in systems of pinned billiard balls. \bigskip

\section{Methodology}
\label{Sec:Methodology}
\subsection{Physics of Collisions}
In classical mechanics, there are three types of collisions: elastic, inelastic, and perfectly inelastic collision. A collision is modeled with respect to the concepts of momentum and energy. Momentum is a product of mass and velocity, measuring the motion of an object as in Definition \ref{def:Momentum}. 
\begin{definition}
    Momentum of an object is defined by
    \begin{align*}
        \mathbf{p} &= \int_{V} \rho \mathbf{v}dV
    \end{align*}
    where $\rho$ is the density of mass and $\mathbf{v}$ is the velocity vector.
    \label{def:Momentum}
\end{definition}
When $\rho$ and $\mathbf{v}$ are constant in Definition \ref{def:Momentum}, momentum simplifies to $\mathbf{p} = m \mathbf{v}$ where $m$ is the total mass. Note that we neglect body forces, i.e. an external field such as gravity as modeled by \cite{BurdzyAndRizzolo2016RandomLorentzGasVariableDensityWithGravity}, in our efforts.

Another quantity that governs the physics of collisions is kinetic energy, described in Definition \ref{def:Kinetic}.
\begin{definition}
    Kinetic energy of an object is defined by
    \begin{align*}
        K =\int_{V} \frac{1}{2} \rho | \mathbf{v}|^2 dV
    \end{align*}
    where $\rho$ is the density of mass and $\mathbf{v}$ is the velocity vector.
    \label{def:Kinetic}
\end{definition}
As in the case of momentum, we can simplify this to $K = \frac{1}{2} m \mathbf{v}^2$ when $\rho$ and $\mathbf{v}$ are uniform under the volume $V$.

\subsection{Collision Dynamics}
\label{Sec:CollisionDynamics}
We emphasize that collisions considered are highly idealized, and hence, more suited to mathematical analysis and certain fundamental physical systems. Collisions are taken as instantaneous, with full momentum transfer. However, as this pinned system enforces contact between all bodies strongly, it behaves as a structural frame system to a certain extent and hence we observe ``load-paths'' in a dynamic form by the trigonometric splitting of these collisions on momentum (which must be conserved by the governing equations).  

While it is of some interest to investigate inelastic, hyper-elastic, and elasto-plastic collisions, i.e. for sticky collisions, large deformation collisions, and collision of bodies that better simulate ivory and plastics used in billiard balls under extreme conditions, this is not undertaken herein.

The collision of objects is elastic if it occurs under the conservation of momentum, energy, and mass (the continuity equation, enforced by our static ball masses). That is, in the two bodies case, these equations must be satisfied:
\begin{align}
    \label{eq:ElasticTwoBody}
    m_1 \mathbf{v}_1 + m_2 \mathbf{v}_2 &= m_1' \mathbf{v'}_1 + m_2' \mathbf{v'}_2 \text{ and }\\
    \notag \frac{1}{2} m_1 \mathbf{v_1}^2 + \frac{1}{2} m_2 \mathbf{v_2}^2 &= \frac{1}{2} m_1' \mathbf{v_1'}^2 + \frac{1}{2} m_2' \mathbf{v_2'}^2 \,\,, 
\end{align}
where $m_i$, $\mathbf{v_i}$ are the mass and velocity of object $i$ before the collision and $m_i'$, $\mathbf{v_i}'$ are the mass and velocity of object $i$ after the collision. We assume that the mass of each object is constant throughout the collision.

Alternative to the above system in Equations \ref{eq:ElasticTwoBody}, we can compute the velocities after the collision following Algorithm \ref{alg:cap}:
\begin{algorithm}[htbp]
    \caption{An Elastic Collision of Two Bodies.}
    \begin{algorithmic}[htbp]
        \Require $m_1, m_2, \mathbf{v}_1, \mathbf{v}_2, \mathbf{d} \,\,,$\\
        $\quad \text{ where } \mathbf{d} \text{ is a unit vector of the collision direction.}$
        \Ensure $\mathbf{v_1'}, \mathbf{v_2'}$
        \State $\mathbf{v_1}^{\parallel} \gets (\mathbf{v_1} \cdot \mathbf{d}) \mathbf{d}$
        \State $\mathbf{v_1}^{\perp} \gets \mathbf{v_1} - \mathbf{v_1}^{\parallel}$
        \State $\mathbf{v_2}^{\parallel} \gets (\mathbf{v_2} \cdot \mathbf{d}) \mathbf{d}$
        \State $\mathbf{v_2}^{\perp} \gets \mathbf{v_2} - \mathbf{v_2}^{\parallel}$
        \\
        \State $\mathbf{v_1'} \gets \mathbf{v_1}^{\perp} + \mathbf{v_2}^{\parallel}$
        \State $\mathbf{v_2'} \gets \mathbf{v_2}^{\perp} + \mathbf{v_1}^{\parallel}$
    \end{algorithmic}
    \label{alg:cap}
\end{algorithm}

Remark: The collision of more than two bodies is not as simple as the case of two.

In real-world physics, when an object has non-zero velocity, its position changes over time. Simulating collisions under this assumption is complex because we need to know each object's exact position to determine if they collide. A collision happens when two objects touch each other.

In the pinned billiard ball system, the balls are fixed in place, but they have a kind of `pseudo-velocity' based on the rules we set. Neighboring balls can `collide' in a way similar to an elastic collision, where this pseudo-velocity can pass from one ball to another.

For this study, we are looking at these collisions in a 2-dimensional flat space. The pinned billiard balls are arranged in a triangular lattice, which is the tightest way to arrange them in this space. The simulation follows Monte Carlo simulations. In each time step, we uniformly and randomly choose a pair of adjacent balls and perform the collision as defined above. 

The simulation is described below:
\begin{algorithm}[htbp]
    \caption{An Elastic Collision of Two Bodies.}
    \begin{algorithmic}[htbp]
        \Require $Ps$ (a list of the balls' positions), $Vs$ (a list of the balls' initial velocities), $Adjs$ (a list of pairs of indices of adjacent balls), $N$ (a number of time steps).
        \State {$n\gets 0$}
        
        \While{$n < N$}
        \State $i \gets \text{Uniformly Sample from 1,2,...,N}$
        \State $j_1, j_2 \gets Adjs_i$ \\ 
        \State $\mathbf{p_1} \gets Ps_{j_1}$
        \State $\mathbf{p_2} \gets Ps_{j_2}$ \\
        \State $\mathbf{v_1} \gets Vs_{j_1}$
        \State $\mathbf{v_2} \gets Vs_{j_2}$ \\
        \State $\mathbf{d} \gets \frac{\mathbf{p_1} - \mathbf{p_2}}{\|\mathbf{p_1} - \mathbf{p_2}\|}$
        \If{$d \cdot (\mathbf{v_1} - \mathbf{v_2}) < 0$}\Comment{check if collision occurs}
            \State $\mathbf{v_1}^{\parallel} \gets (\mathbf{v_1} \cdot \mathbf{d}) \mathbf{d}$
            \State $\mathbf{v_1}^{\perp} \gets \mathbf{v_1} - \mathbf{v_1}^{\parallel}$
            \State $\mathbf{v_2}^{\parallel} \gets (\mathbf{v_2} \cdot \mathbf{d}) \mathbf{d}$
            \State $\mathbf{v_2}^{\perp} \gets \mathbf{v_2} - \mathbf{v_2}^{\parallel}$
            \\
            \State $\mathbf{Vs_{j_1}} \gets \mathbf{v_1}^{\perp} + \mathbf{v_2}^{\parallel}$
            \State $\mathbf{Vs_{j_2}} \gets \mathbf{v_2}^{\perp} + \mathbf{v_1}^{\parallel}$
            \EndIf
            \EndWhile
    \end{algorithmic}
    \label{alg:cap_adjs}
\end{algorithm}

For a billiard ball $b$ with mass $m_b$ and velocity vector $\mathbf{v}_b$, we may take its velocity magnitude with Euclidean norm operator $|| \circ ||$ where $\circ$ is a placeholder for any input vector quantity (e.g. velocity). The kinetic energy of the full, discretely represented billiard ball system is then
\begin{align*}
    \mathcal{E}_{\text{Kinetic}} = \frac{1}{2} \sum_b{m_b ||\mathbf{v}_b||^2}\,\,,
\end{align*}
when viewed in the classical sense. We take each ball's mass to be a scalar value of 1, stated as $m_b = 1 \quad \forall b$ and the initial total energy to be one hundred. Note that this is similar to ``lumped mass" formulations seen in computational mechanics; there is potential to expand to full-mass matrix representations for more nuanced system behavior. 

A respectable numerical solver for pinned billiard ball problems is expected to maintain kinetic energy, $\mathcal{E}_{\text{Kinetic}}$, throughout a simulation's duration unless damping is introduced, absorption is present at boundaries, an instability is present, or floating-point arithmetic errors in the underlying code are non-negligible. In our Python implementation, all computation is done in double-precision, a reasonable avenue for a stiff numerical problem such as pinned billiard systems, though we note that higher precision implementations may be pursued and that underlying precision optimization could be undertaken. Further, any future computationally parallel implementation efforts may benefit from exploring single- and mixed-precision implementation to attain accelerated performance on Graphics Processing Units (GPUs). 
\subsection{Half-Space Configuration}
\label{Sec:HalfSpaceConfig}
In Euclidean space, a half-space is defined to be a part of space constructed from dividing an affine space by a hyperplane.
\begin{definition}
    In n-dimensional Euclidean Space, a closed-half space is defined to be a set of all $\mathbf{x} \in \mathbb{R}^n$ such that $\mathbf{a}^T\mathbf{x} \leq b$, where $\mathbf{a} \in \mathbb{R}^n$ and $b \in \mathbb{R}$.
\end{definition}
In 2-dimensional space, a half-space is a half-plane. Additionally, we define a wall or a boundary of a half-space to be a set in the form $\{x \in \mathbb{R}^n: \mathbf{a}^T\mathbf{x} = b\}$

In this project, we consider a system where pinned billiard balls are tightly placed on a 2-dimensional half-space with the initial velocity of all balls to be $\mathbf{0}$ except for one ball on the boundary whose velocity is non-zero and directed toward the balls.

\begin{figure}[htbp]
    \centering
    \includegraphics[width=0.4\columnwidth, trim = {0cm 8cm 8.125cm 8cm}, clip, angle=90]{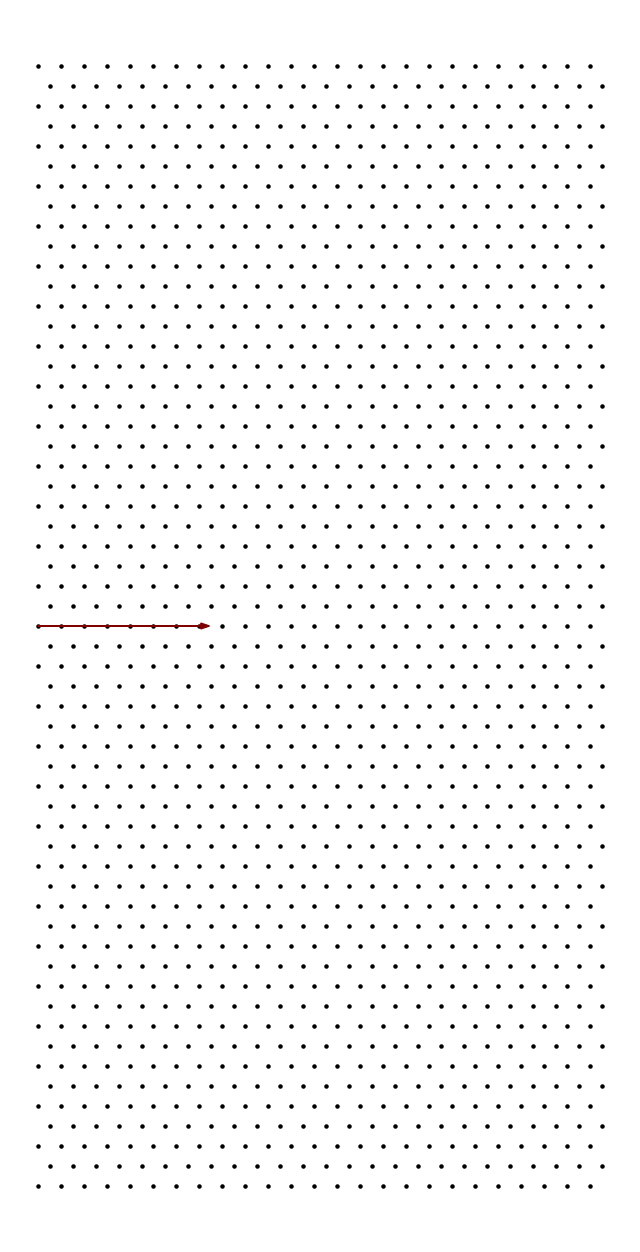}
    \captionsetup{font=small}
    \caption{An approximated 2-dimensional half-space with initial velocities, illustrating the configuration of pinned billiard balls.}
    \label{fig:half_space}
\end{figure}

Because the half-plane is infinite, it is implausible to construct a full lattice on a computer. So, we approximate the half-space using rectangular configurations with the following dimensions: $10\times 5$, $60\times 30$, $100\times 50$, $200\times 100$, $300\times 150$, $400\times 200$, and $500\times 250$. We will display visualizations from the $100\times 50$ case. In this system, the rectangle contains 1425 balls. The approximated half-space and initial velocities are shown in Figure \ref{fig:half_space} with some vertical extents clipped for visual brevity.

Under the laws of conservation of energy and linear momentum \ref{def:Momentum}, when two objects collide with one object of infinite mass, the momentum of the colliding object would direct in the opposite of the colliding direction. Therefore, we conjecture that after a number of collisions, a large quantity of velocities will direct outward from the boundary.
    
\subsection{Flat Torus Configuration}
\label{Sec:TorusConfig}
\begin{figure}[htbp]
    \centering
    \includegraphics[scale = 0.2]{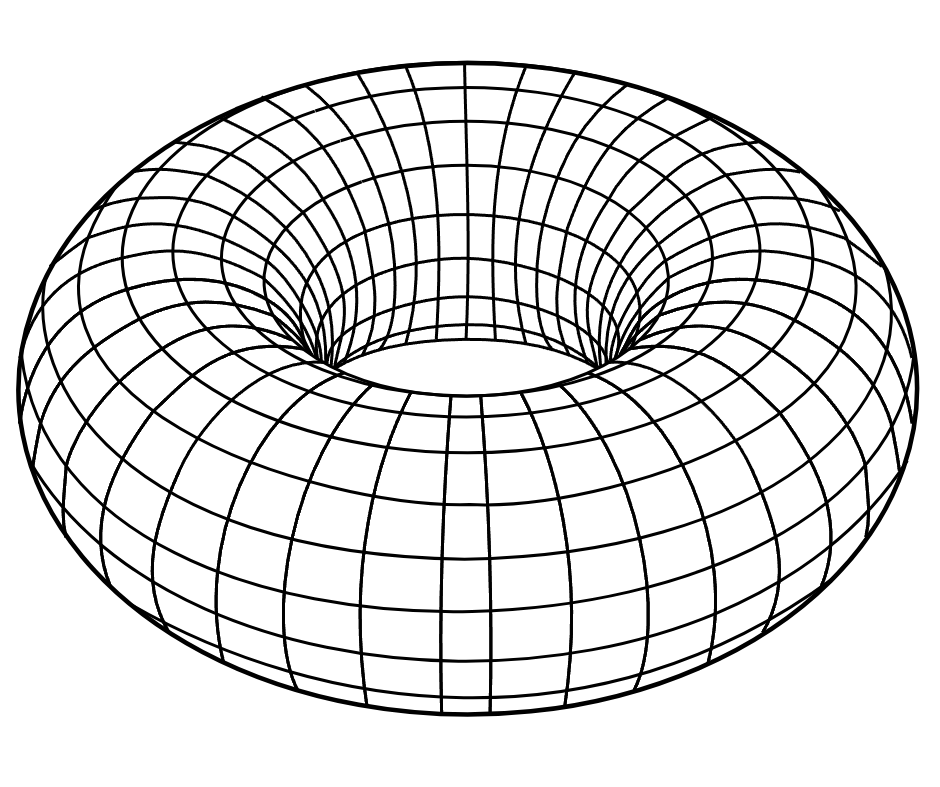}
    \captionsetup{font=small}
    \caption{A ring torus illustration. Its spatial domain "wraps-around" itself, producing periodic boundaries when presented in a flattened configuration.}
    \label{fig:ring_torus}
\end{figure}

A torus is defined be a surface generated by rotating a circle in 3-dimensional Euclidean space, shown in Figure \ref{fig:ring_torus}. A ring torus is the surface whose axis of revolution does not touch the circle.


Consider that the ring torus can be unfolded into a bounded 2-dimensional plane so-called flat torus. In spite of being bounded, the flat torus has a periodic boundary, following from its construction from the ring torus.
We consider a system where pinned billiard balls are placed on a flat torus. The periodic boundary allows balls from opposite sides of the plane to collide. 
\begin{figure}[htbp]
    \centering
    \includegraphics[width=0.55\columnwidth, trim = {0cm 1cm 0cm 0cm}, clip]{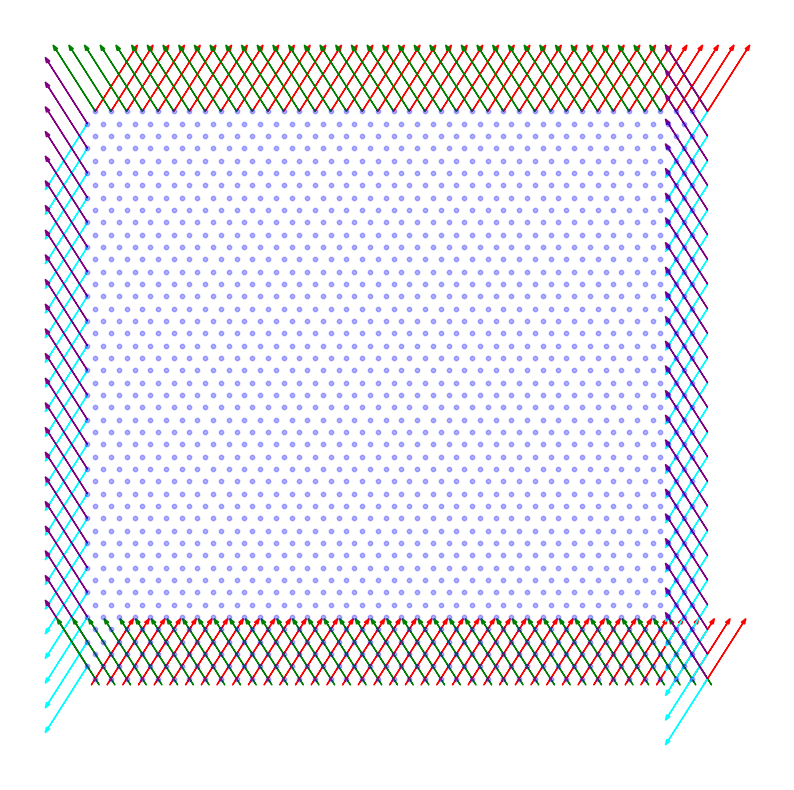}
    \captionsetup{font=small}
    \caption{Periodic boundaries applied on the borders of a flattened torus domain.}
    \label{fig:torus}
\end{figure}

Since the boundary apparently disappears, unlike in the typical case where balls on the boundary have fewer neighbors, each ball can collide in all directions. As such, the system may experience an infinite number of collisions. In addition, as the pseudo-velocities circulate, they tend to be uniformly distributed.

Theorem \ref{th:uni_dist_sphere}, when $n = 2$ in 2-dimensional space,  the x-components and y-components of the velocity vectors are independent and identically distributed random variables, standard normal. 
\begin{theorem}
    Let $(X_1, X_2, \ldots, X_n)$ be uniformly distributed on an $n$-dimensional sphere $\mathbb{S}^{n-1}$. Then, $\sqrt{n}X_1$ converges to standard normal distribution as $n\to \infty$.
    \label{th:uni_dist_sphere}
\end{theorem}

\begin{figure}[htbp]
    \centering
    \includegraphics[width=0.5\columnwidth, trim = {0.5cm 13.5cm 13cm 1cm}, clip, angle=0]{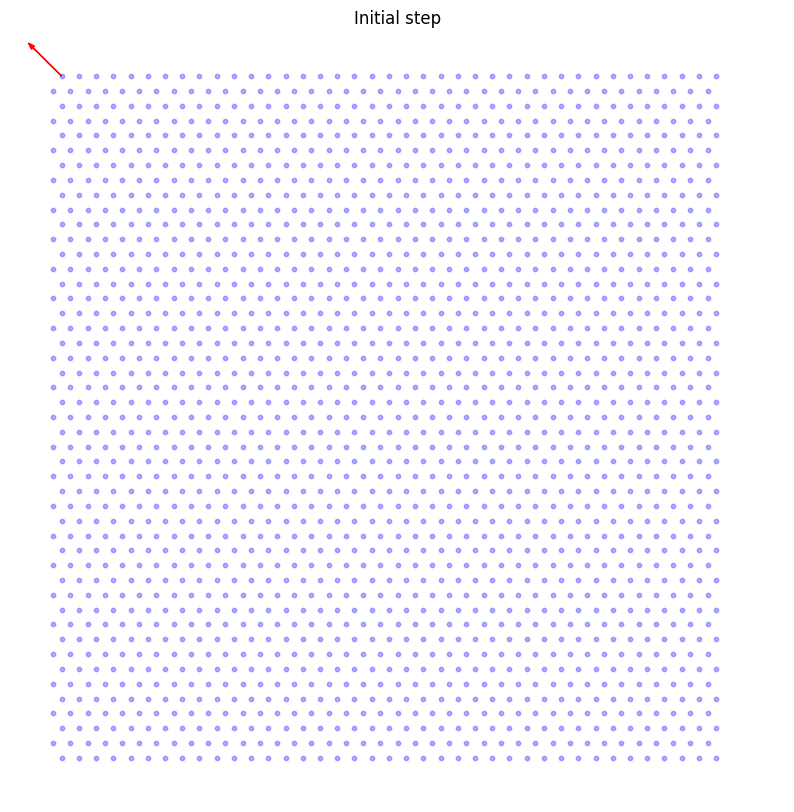}
    \captionsetup{font=small}
    \caption{An approximated flattened-torus pinned billiard ball system with an initial velocity on a single ball at the upper-right corner. Balls are in a staggered fabric and are always in contact with neighbors.}
    \label{fig:init_step_k}
\end{figure}

For a flat torus system, we undertake a simulation in our open-source Python code-base. For the initial state, we select one billiard ball within an array of billiard balls to possess an initial velocity vector, while all other balls are initially stationary. This initial state is visualized in Figure \ref{fig:init_step_k}, with the velocity vector visualized in red and the billiard balls represented as discrete dots. Note that despite not visually touching one-another, the billiard balls are initially in contact with all their adjacent neighbors and remain so throughout the simulation.

Then, we numerically compute the distributions of x- and y- components of the velocity vectors. Furthermore, for pairs of adjacent balls, we compute the correlation of the velocities projected on a line passing through the centers of the pairs.

With the applied methodologies defined, we now continue into the numerical results produced and contributions resolved from these efforts in the study of pinned billiard ball systems. 
   
\section{Numerical Results}
\label{Sec:NumericalResults}
Numerical results for pinned billiard ball systems, spanning hundreds of simulations with applied methodologies noted in Section \ref{Sec:Methodology}, are discussed herein. Results are segmented between systems within a \textbf{(i)} half-space and a \textbf{(ii)} flattened torus-space. We demonstrate their respective tendencies to accumulate energy near enforced boundaries and to distribute collision and velocity components along normal Gaussian distributions over appropriate time-scales.

\subsection{Half-plane Configuration}
\label{Sec:NumericalResults_HalfSpace}
\begin{figure}[htbp]
    \centering
    \includegraphics[width=0.4\columnwidth, trim = {0cm 0cm 0cm 0cm}, clip, angle=90]{halfspac.png}
    \vfill
    \includegraphics[width=0.4\columnwidth, trim = {0cm 0.5cm 0cm 1cm}, clip, angle=90]{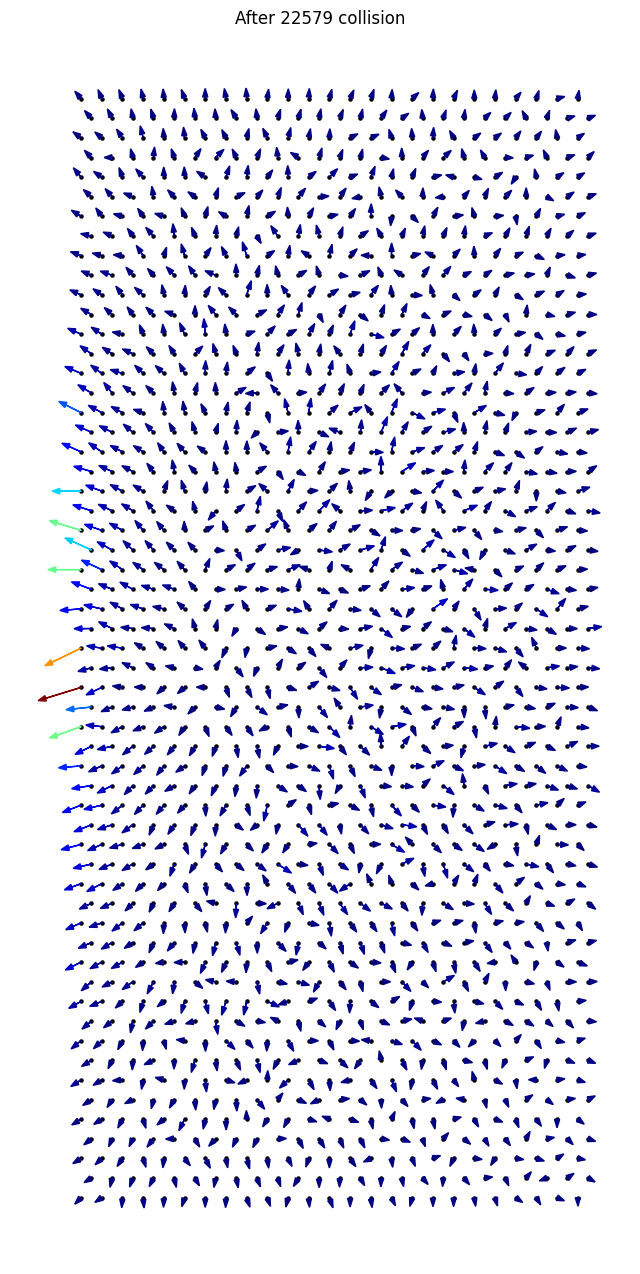}
    \vfill 
    \includegraphics[width=0.4\columnwidth, trim = {0cm 0.5cm 0cm 1cm}, clip, angle=90]{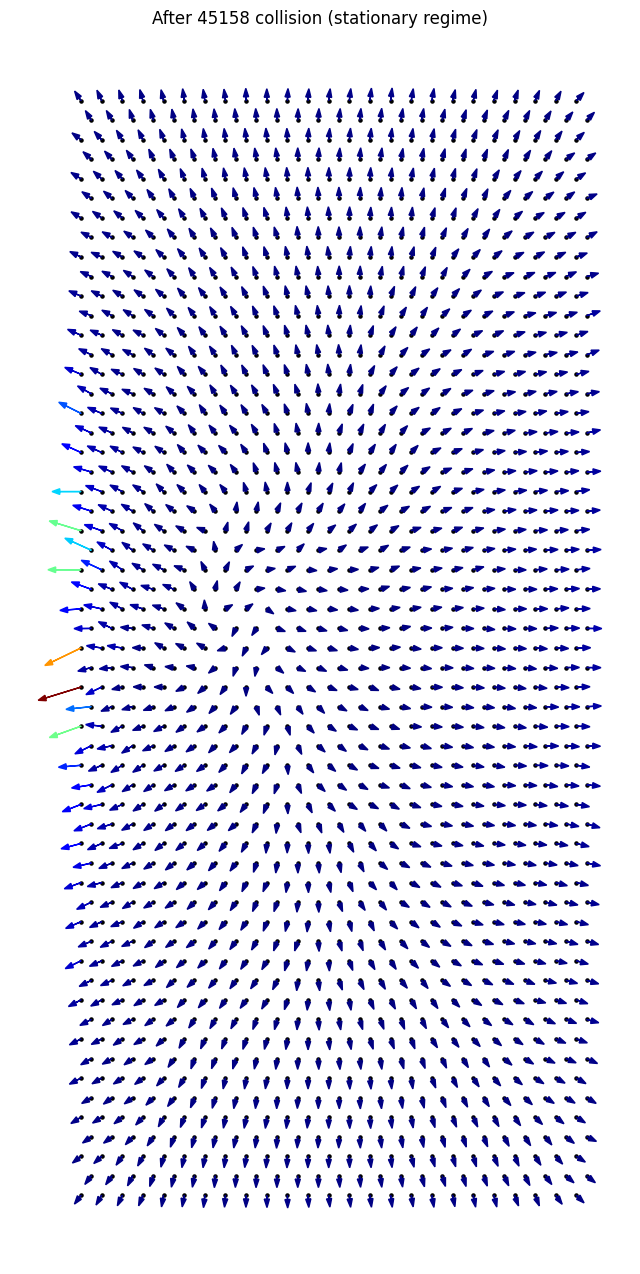}
    \captionsetup{font=small}
    \caption{States of half-plane pinned billiard ball simulation following set numbers of collisions. The $x$-component is up, $y$-component is to the left. Velocity (visualized as colored vectors) shown per ball (visualized as grey dots). All billiard balls remain in contact throughout the simulations. (Top) Initial state. (Middle) State after 211,417 collisions. (Bottom) State after 422,834 collisions.}
    \label{fig:num_collisions_half_2Panel}
\end{figure}

Figure \ref{fig:num_collisions_half_2Panel} shows arrangements of the billiard balls in an ordered fabric occupying a half-plane, where velocities are visualized at simulation snapshots when the collision count reaches 211,417 collisions and 422,834 collisions.

\begin{figure}[htbp]
    \centering
        \includegraphics[width=0.9\columnwidth, trim = {0cm 0cm 0cm 0cm}, clip]{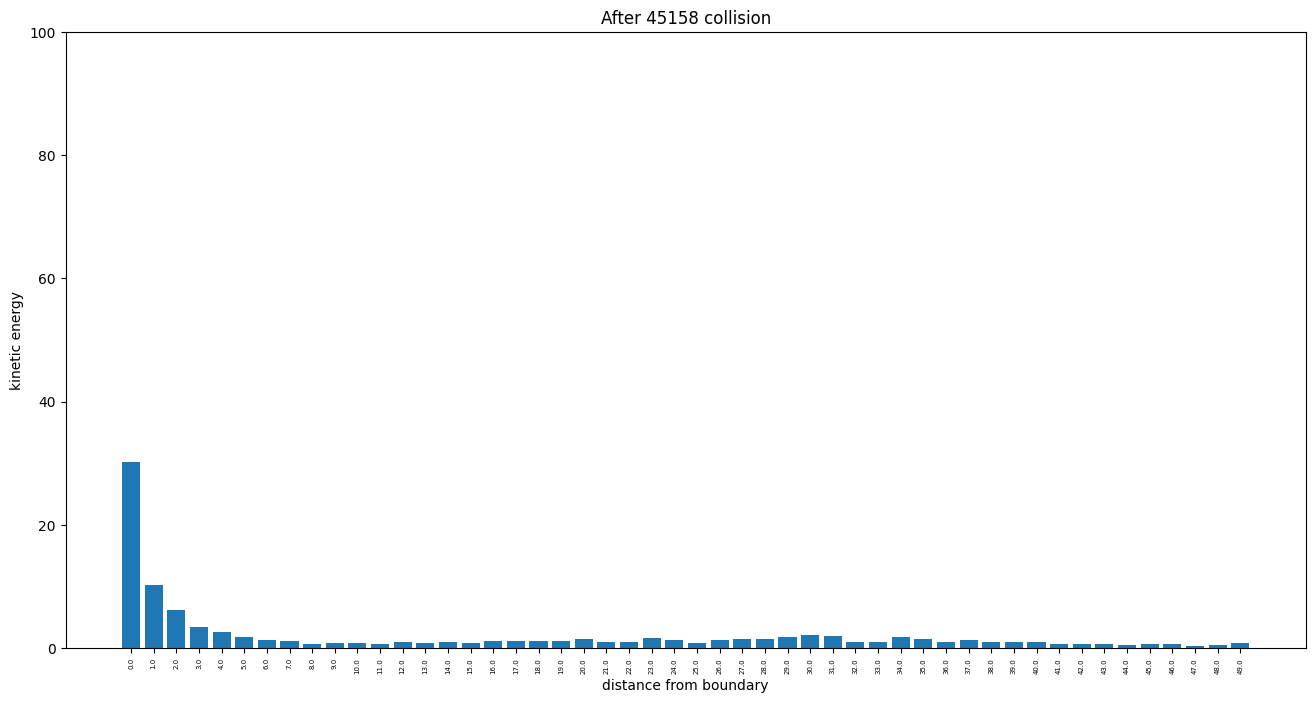}
        \centering \includegraphics[width=0.9\columnwidth, trim = {0cm 0cm 0cm 0cm}, clip]{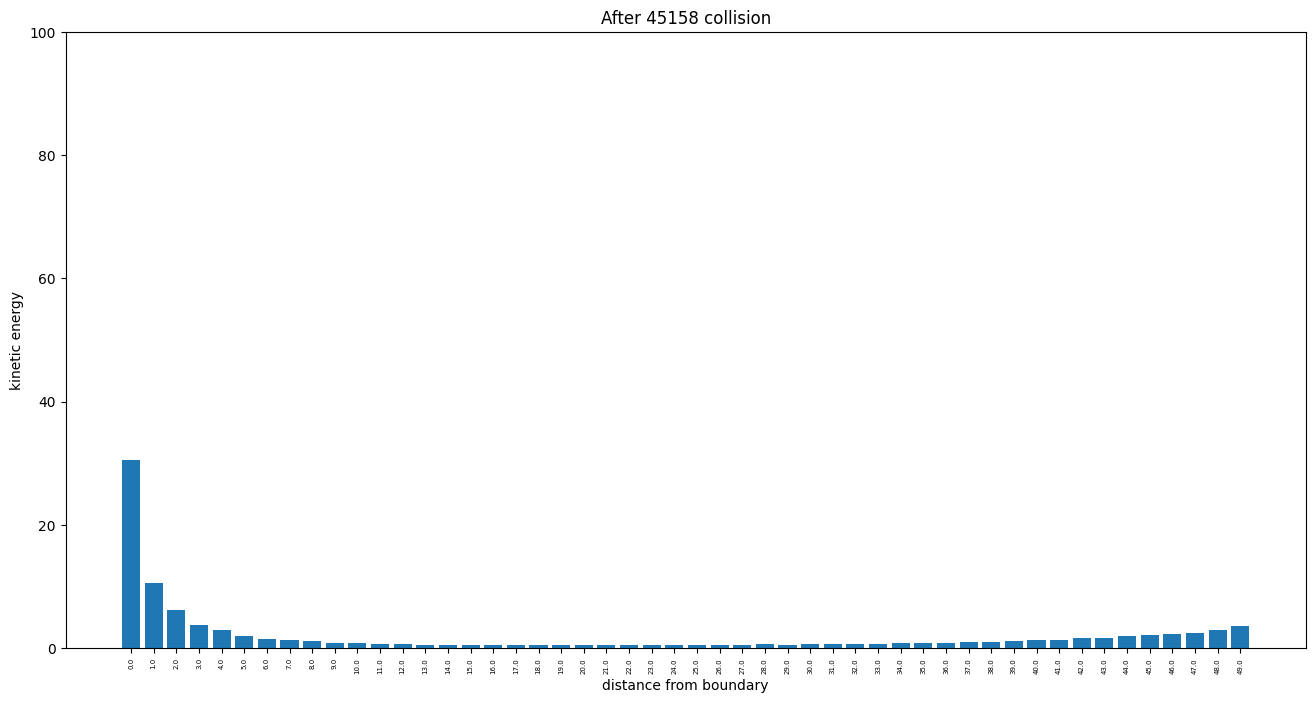}
    \captionsetup{font=small}
    \centering   
    \caption{The kinetic energy of the system at the distance away from the boundary after 22,579 collisions (Top) and in the stationary regime following 45,158 collisions (Bottom).}
    \label{fig:Bar}
\end{figure}

The bar plots (Figure \ref{fig:Bar}) display the kinetic energy of balls at increasing distances from the boundary after 22,579 and 45,158 collisions.

The result from numerical simulation aligns with a conjecture that a majority of velocity will cluster around the boundary in the outward direction. Approximately 30 percent of the energy in the system was concentrated at the boundary. In addition, a small amount of the velocity travels inward and settles at another boundary of the system.

\begin{figure}[htbp]
    \centering
        \includegraphics[width=0.95\columnwidth, trim = {0cm 0cm 0cm 0cm}, clip]{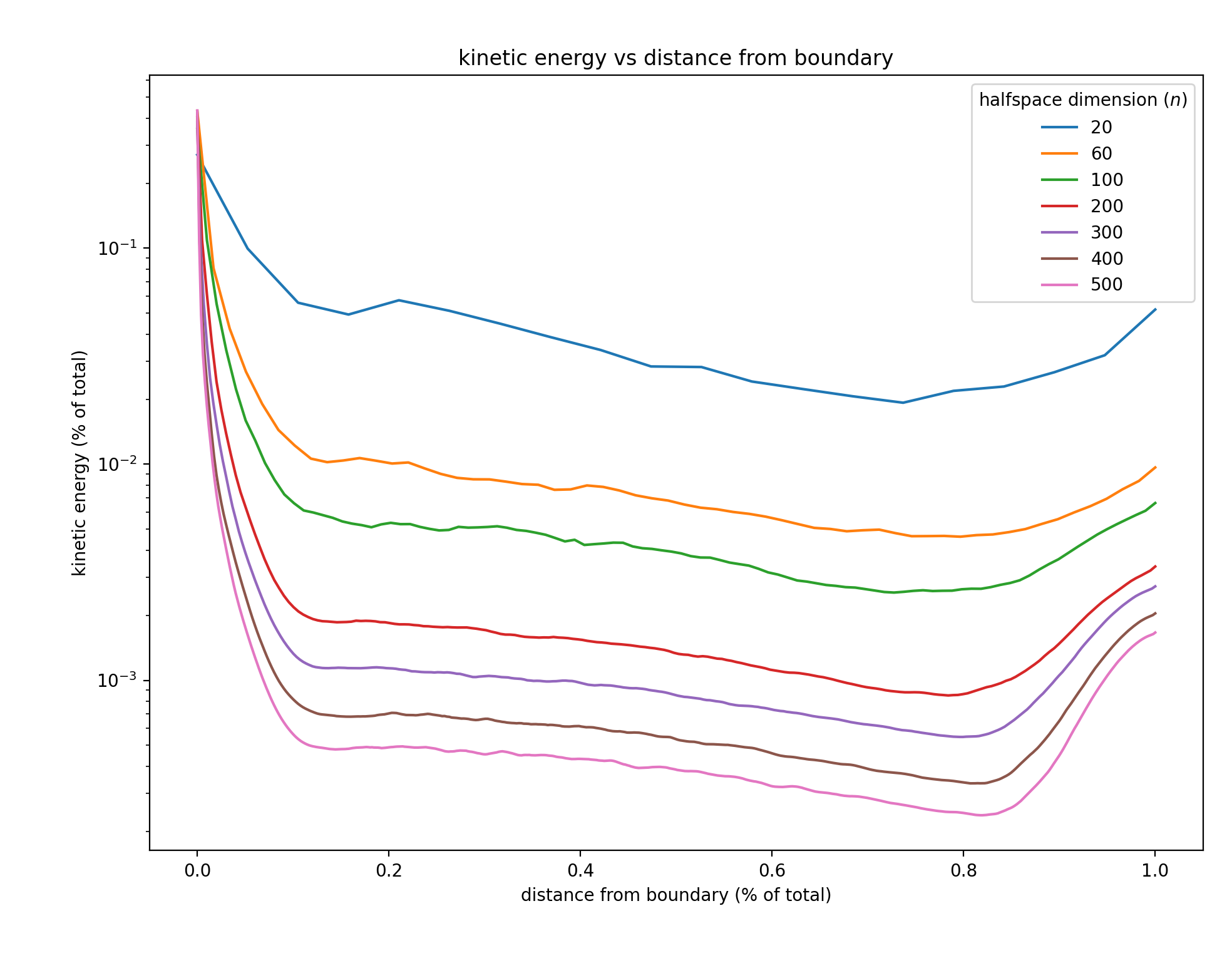}
    \captionsetup{font=small}
    \centering   
    \caption{The distribution of kinetic energy in each system as distance from the boundary increases. From top to bottom, the lines represent the respective cases: $10\times 5$, $60\times 30$, $100\times 50$, $200\times 100$, $300\times 150$, $400\times 200$, and $500\times 250$.}
    \label{fig:Colins}
\end{figure}

Furthermore, we ran simulations for various sizes of the systems; the distributions of energies for each system are shown in Figure \ref{fig:Colins}.

Similar to the previous bar plot, most of the energy is concentrated around the boundary. However, as the system gets bigger, the energies far away from the boundary get smaller. This aligns with the conjecture that the energy will cluster only on the boundary for an infinite system.

Numerical results garnered in our code-base for pinned billiard balls in a half-space with a boundary align with a conjecture that claims that the majority of velocity vectors (each individual ball, $b$, possessing one vector) will cluster around boundaries (i.e. system constraints) with an outward bearing. This suggests that either there is a natural tendency for stiff systems to radiate energy outwards, or to achieve a balanced, steady system with the assistance of immovable boundaries at the billiard ball-rigid boundary interface. Approximately 30 percent of the energy in the system is seen to accumulate here for the time-intervals evaluated and billiard ball configurations tested. In addition, a small amount of the velocity travels inward and settles at another end of the boundary, though this may be a transient property that will fade given appropriate time. 

\subsection{Torus Configuration}
\label{Sec:NumericalResults_Torus}
In the flat torus system, whose configuration is described in Section \ref{Sec:TorusConfig}, our numerical results suggest the system reaches a steady-state after 45,158 collisions. 


\begin{figure}[htbp]
    \centering
    \includegraphics[width=0.6\columnwidth, trim = {1cm 1cm 1cm 1cm}, clip, angle=0]{initial_step_k.png} \\ 
    \includegraphics[width=0.6\columnwidth, trim = {1cm 1cm 1cm 1cm}, clip, angle=0]{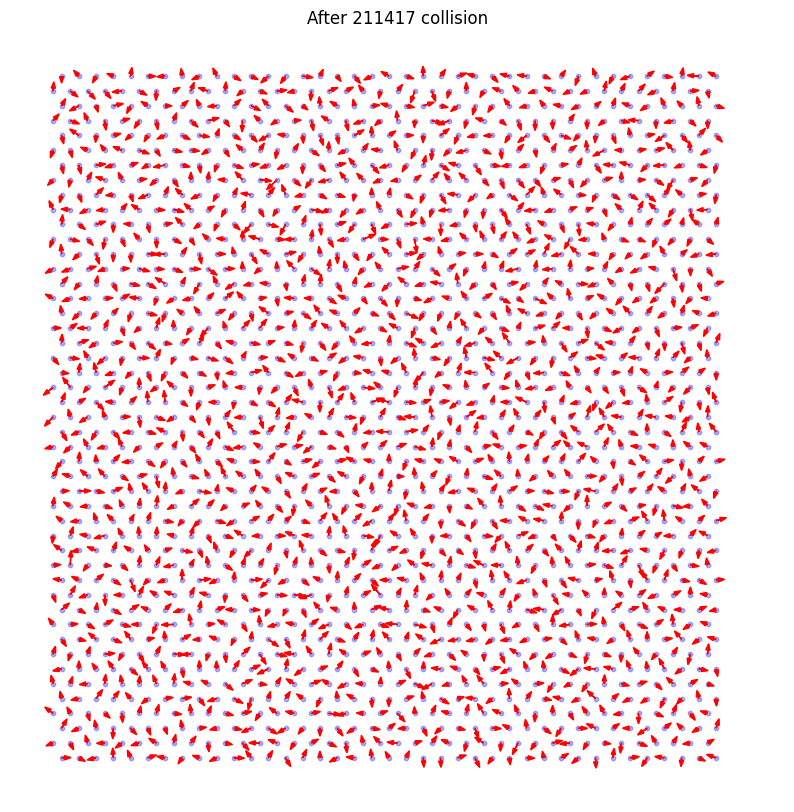}\\
    \includegraphics[width=0.6\columnwidth, trim = {1cm 1cm 1cm 1cm}, clip, angle=0]{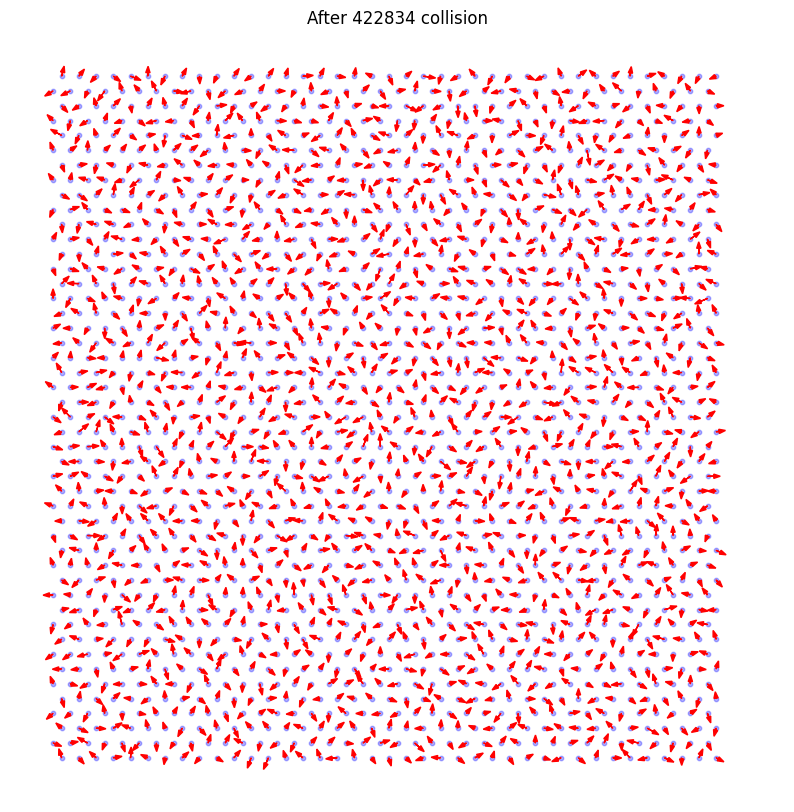}
    \captionsetup{font=small}
    \caption{States of the torus configuration simulation following set numbers of collisions. Velocity (visualized as red vectors) shown per ball (visualized as grey dots). All billiard balls remain in contact throughout the simulations. (Top) Initial state. (Middle) State after 211,417 collisions. (Bottom) State after 422,834 collisions.}
    \label{fig:num_collisions_torus_2Panel}
\end{figure}

\begin{figure}[htbp]
    \centering
    \includegraphics[width=0.9\columnwidth, trim = {0cm 0cm 0cm 0.65cm}, clip]{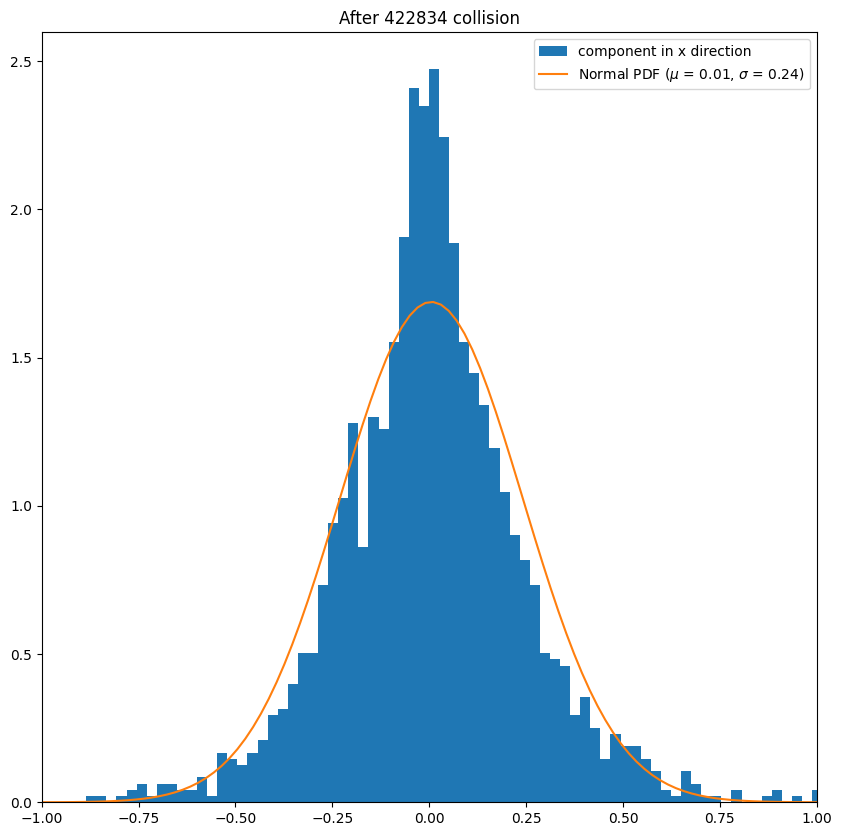}
    \includegraphics[width=0.9\columnwidth, trim = {0cm 0cm 0cm 0.65cm}, clip]{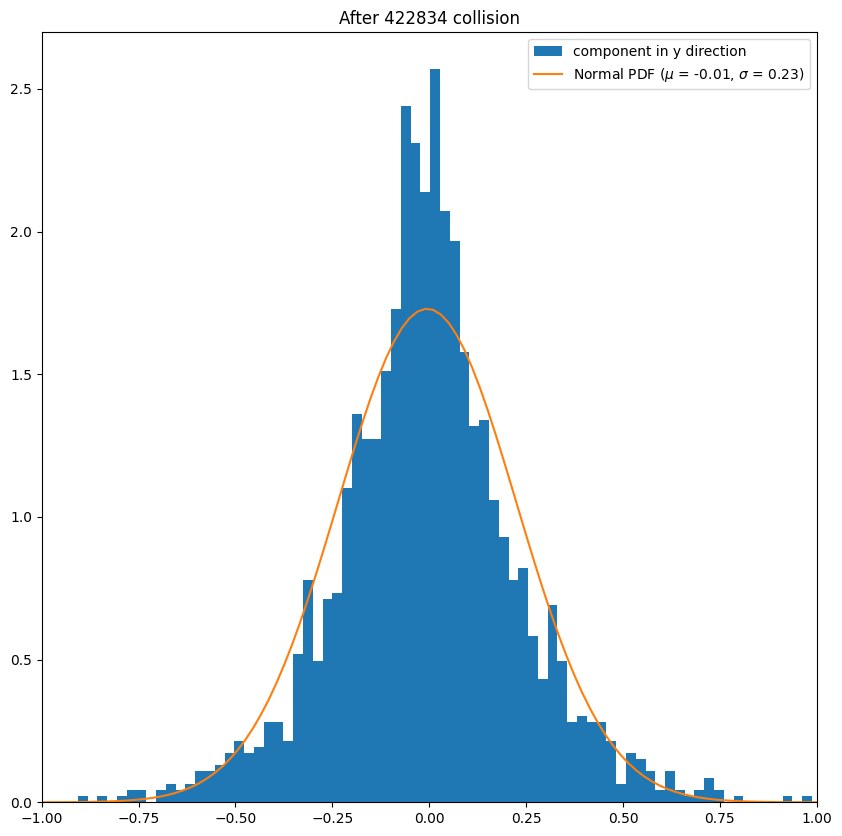}
    \captionsetup{font=small}
    \centering
    \caption{Histogram distributions of velocity components in the torus configuration following set numbers of collisions (211,417). Empirical data binned as a histogram with bin-widths of 0.025 over the normalized -1 to 1 span. A fitted Gaussian distribution is super-imposed in orange. All billiard balls remain in contact throughout the  simulations. (Left) X-component results. (Right) Y-component results.}
    \label{fig:histogram_2k_2Panel}
\end{figure}

\begin{figure}[htbp]
    \centering
    \includegraphics[width=0.9\columnwidth, trim = {0cm 0cm 0cm 0.65cm}, clip]{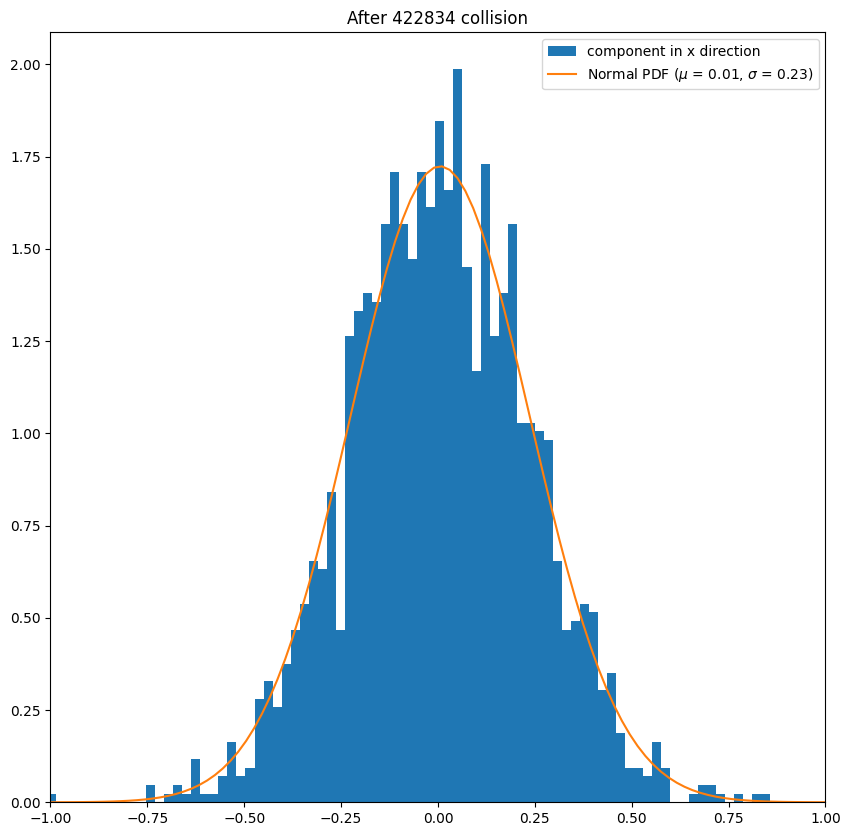}

    \includegraphics[width=0.9\columnwidth, trim = {0cm 0cm 0cm 0.65cm}, clip]{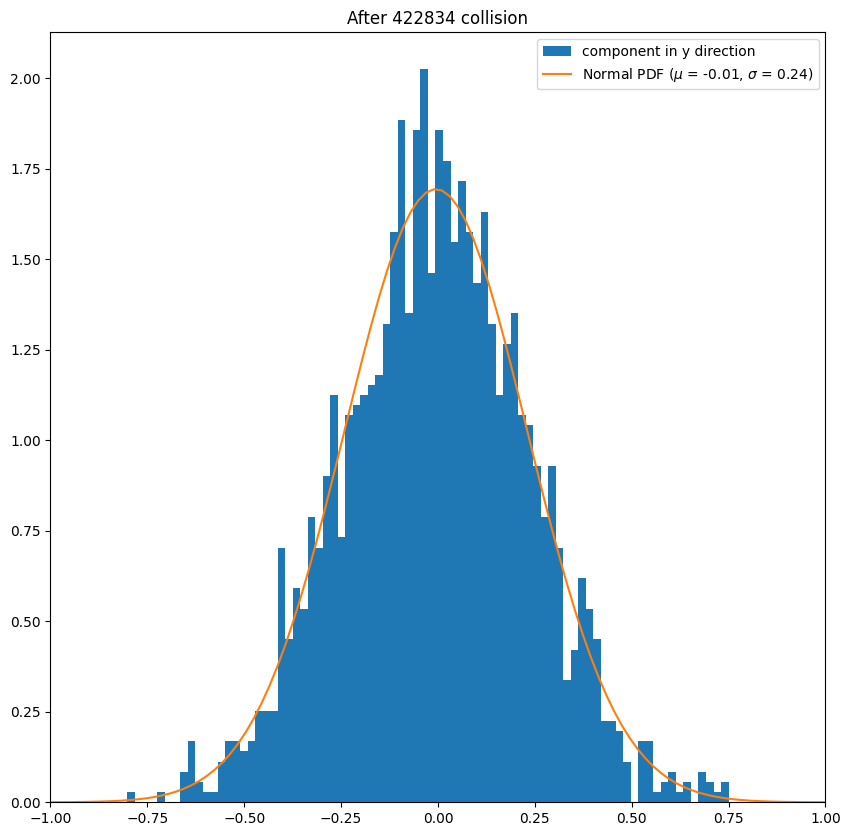}
    \captionsetup{font=small}
    \centering
    \caption{Histogram distributions of velocity components in the torus configuration following set numbers of collisions (422,834). Empirical data binned as a histogram with bin-widths of 0.025 over the normalized -1 to 1 span. A fitted Gaussian distribution is super-imposed in orange. All billiard balls remain in contact throughout the simulations. (Left) X-component results. (Right) Y-component results.}
    \label{fig:histogram_4k_2Panel}
\end{figure}

The distributions of x- and y-components of the velocity vectors gradually become normal over appropriate time-scales in the simulation, shown in Figures \ref{fig:histogram_2k_2Panel} and \ref{fig:histogram_4k_2Panel}, suggesting that initially non-normal perturbed behavior is transient for pinned-billiard ball systems and will tend to dissipate. The mean of the distribution is zero, by symmetry, i.e. $\mu_{\mathrm{Mean}} = 0$.

\begin{figure}[htbp]
    \centering
    \includegraphics[width=0.8\columnwidth, trim = {0cm 0cm 0cm 0.75cm}, clip]{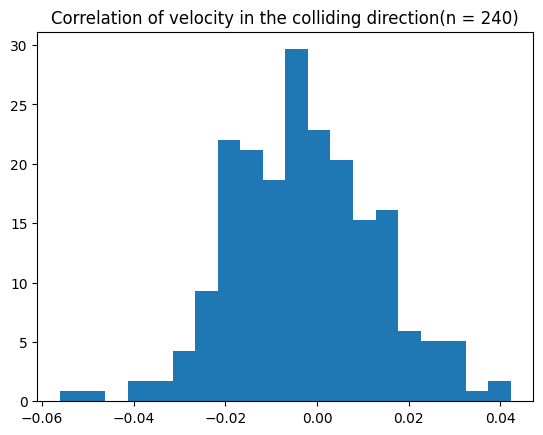}
    \captionsetup{font=small}
    \caption{Distribution of correlations of collisional velocity components calculated from 240 simulations. Empirical data binned as a histogram with bin-widths of 0.02 over the typical -1 to 1 range.}
    \label{fig:correlation_collision}
\end{figure}

\begin{figure}[htbp]
    \centering
    \includegraphics[width=1.05\columnwidth, trim = {0cm 0cm 0cm 0.75cm}, clip]{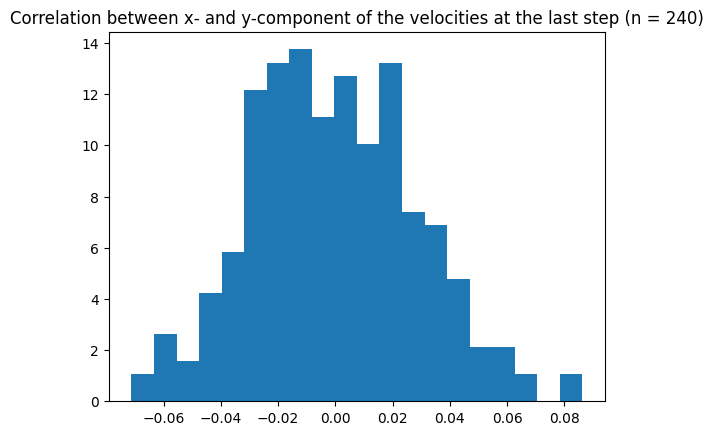}
    \captionsetup{font=small}
    \caption{Distribution of correlations of x-y velocity components calculated from 240 simulations. Empirical data binned as a histogram with bin-widths of 0.02 over the typical -1 to 1 range.}
    \label{fig:correlation_xy}
\end{figure}

Additionally, we ran the simulation multiple times to calculate the correlations between the collision's components (Figure \ref{fig:correlation_collision}) and x-y components (Figure \ref{fig:correlation_xy}). We found that both correlation values were close to zero. 

\section{Discussion}
\label{Sec:Discussion}
We now consider the significance of our results for the two infinite configurations. 

As the size of our half-plane model increases, energy becomes more concentrated toward the boundary. So, we suspect that in the actual half-plane system, we would see almost all of the energy in the boundary itself. Furthermore, the velocities primarily pointed toward the boundary, indicating that the boundary is the outlet of energy in the system. These results match those from finite configurations, despite there being no true terminal state of the half-plane configuration.  

In the torus configuration, we found that the x and y components of the velocities were normally distributed. This result, combined with the near-zero correlations between the x-y components and collisional components, suggests that all components of all velocities are independent normal in the stationary regime. This result supports the physics postulate of equi-distribution of energy. 
    
\section{Conclusions}
\label{Sec:Conclusions}
In this paper, we investigate a system of pinned billiard balls, which is a simplified version of collision systems. In this system, every collision object is fixed (i.e. pinned) spatially; yet, the velocity of each object can circulate under the governing rules of elastic collisions. We numerically studied two systems, one occupying a \textbf{(i)} half-space and the other occupying a \textbf{(ii)} flat torus.

Arranging balls on a half plane with a sole initial collision on the boundary, we found that, in the stationary regime, a large amount of velocity settled around and pointed outward from the boundary. The numerical results substantiate the conjecture we made regarding this system.

In a flat torus, we study the system whose balls are allowed to collide with other balls on the opposite boundary, which is a type of periodic boundary. We found that after a large number of collisions, the distributions of x- and y-components are normally distributed.

\section{Future Work}
\label{Sec:Future}
Investigations into half-space and flattened torus domains have yielded intriguing, but preliminary, results. There is a vested interest in applying our pinned billiard formulation to three primary examples of curved spaces, namely \textbf{(i)} flat-space, \textbf{(ii)} spherical space, and \textbf{(iii)} hyperbolic-space. While the former flat-space is the most obviously useful for applied cases and simplified analysis, both spherical and hyperbolic spaces should not be discounted. There is an active discourse in the field of physics regarding whether or not our universe is best approximated, or truly is, characterized by one of the three spaces. Though viewing a given local regime closely, i.e. a gauge perspective, converges towards a flat-space approximation regardless of which is the overarching behavior of our universe. As particle systems, Brownian motion, etc. can be, to some extent, modeled with a pinned billiard system over some time intervals, it would be of great value to extend our work to spherical and hyperbolic spaces to perhaps shed light on fundamental phenomena in our reality. 

A high-performance parallel implementation is another strong step to take forward. Though our software is open-source and more than adequate for the preliminary study performed herein, this is a clearly parallel problem and thus it should be implemented within a parallel computation architecture to significantly accelerate and grow the simulated problem sets. This may allow for billions of collisions over millions of bodies if well-optimized, producing far stronger conclusions on the potential steady-state behavior of the governing equations.

\section{Acknowledgements}
\label{Sec:Acknowledgements}
We thank the Washington Experimental Mathematics Lab (WXML), University of Washington, Seattle, for providing opportunity and support in this research endeavor. 



\begin{thebibliography}{10}

\bibitem{AthreyaBurdzyDuarte2018PinnedBilliardFoldings}
Jayadev~S. Athreya, Krzysztof Burdzy, and Mauricio Duarte.
\newblock On pinned billiard balls and foldings.
\newblock {\em Indiana Univ. Math. J.}, 72(3):897--925, 2023.

\bibitem{BuragoFerlegerKononenko1998GeometricApproachSemiDispersingBilliards}
D.~Burago, S.~Ferleger, and A.~Kononenko.
\newblock A geometric approach to semi-dispersing billiards.
\newblock {\em Ergodic Theory Dynam. Systems}, 18(2):303--319, 1998.

\bibitem{BuragoFerlegerKononenko1998GlobalBoundNumCollisionsSemiDispersingBilliards}
D.~Burago, S.~Ferleger, and A.~Kononenko.
\newblock Unfoldings and global bounds on the number of collisions for generalized semi-dispersing billiards.
\newblock {\em Asian J. Math.}, 2(1):141--152, 1998.

\bibitem{BuragoFerlegerKononenko1998EstimateOnNumCollisionsSemiDispersingBilliards}
D.~Burago, S.~Ferleger, and A.~Kononenko.
\newblock Uniform estimates on the number of collisions in semi-dispersing billiards.
\newblock {\em Ann. of Math. (2)}, 147(3):695--708, 1998.

\bibitem{BuDu2022}
Krzysztof Burdzy and Mauricio Duarte.
\newblock Upper bound on the number of collisions of pinned billiard balls.
\newblock 2022.
\newblock Arxiv:2203.09013.

\bibitem{BurdzyAndRizzolo2016RandomLorentzGasVariableDensityWithGravity}
Krzysztof Burdzy and Douglas Rizzolo.
\newblock A random flight process associated to a {L}orentz gas with variable density in a gravitational field.
\newblock {\em Stochastic Process. Appl.}, 128(1):79--107, 2018.

\bibitem{BurdzyAndWhite2008MarkovWithProductFormStationaryDistributions}
Krzysztof Burdzy and David White.
\newblock Markov processes with product-form stationary distribution.
\newblock {\em Electron. Commun. Probab.}, 13:614--627, 2008.

\bibitem{Galperin1981MovingLocalRepellingParticles}
G.~A. Galperin.
\newblock Systems of locally interacting and repelling particles that are moving in space.
\newblock {\em Trudy Moskov. Mat. Obshch.}, 43:142--196, 1981.

\bibitem{Galperin1983OnMovingLocalRepellingParticles}
G.~A. Galperin.
\newblock On systems of locally interacting and repelling particles moving in space.
\newblock {\em Trans. Moscow Math. Soc.}, pages 159--214, 1983.

\bibitem{Goebel2018BallsInBanachSpaces}
Kazimierz Goebel and Stanis{\l}aw Prus.
\newblock {\em Elements of geometry of balls in {B}anach spaces}.
\newblock Oxford University Press, Oxford, 2018.

\bibitem{KozlovAndTreshchev1991BilliardsGenetic}
Valeri{\u\i}~V. Kozlov and Dmitri{\u\i}~V. Treshch{\"e}v.
\newblock {\em Billiards: A genetic introduction to the dynamics of systems with impacts}, volume~89 of {\em Translations of Mathematical Monographs}.
\newblock American Mathematical Society, Providence, RI, 1991.
\newblock Translated from Russian by J. R. Schulenberger.

\bibitem{Tabachnikov2005GeometryAndBilliards}
Serge Tabachnikov.
\newblock {\em Geometry and billiards}, volume~30 of {\em Student Mathematical Library}.
\newblock American Mathematical Society, Providence, RI and Mathematics Advanced Study Semesters, University Park, PA, 2005.

\end{thebibliography}

\def\cprime{$'$}

\end{document}